\documentclass[preprint,12pt]{elsarticle}

\usepackage{amssymb,amsthm}
\usepackage{graphicx,amsmath,bm,color,geometry,subfig}
\usepackage{algorithmic,algorithm,diagbox} 
 \usepackage[colorlinks,
             linkcolor=blue,
             anchorcolor=blue,
             citecolor=blue]{hyperref}

\newtheorem{definition}{Definition}

\journal{Journal of Computational Physics}

\begin{document}

\begin{frontmatter}



\title{Koopman operator learning using invertible neural networks\tnoteref{label-title}}
\tnotetext[label-title]{This work is partially supported by the National Natural Science Foundation of China (NSFC) under grant number 12101407, the Chongqing Entrepreneurship and Innovation Program for Returned Overseas Scholars under grant number CX2023068, and the Fundamental Research Funds for the Central Universities under grant number 2023CDJXY-042.}

\author[label-addr1]{Yuhuang Meng\fnref{label-cofirst}}
\ead{mengyh@shanghaitech.edu.cn}

\author[label-addr1]{Jianguo Huang\fnref{label-cofirst}}
\ead{huangjg@shanghaitech.edu.cn}
\fntext[label-cofirst]{Yuhuang Meng and Jianguo Huang contributed equally to this paper.}

\author[label-addr2,label-addr3]{Yue Qiu\corref{label-cor}}
\ead{qiuyue@cqu.edu.cn}
\cortext[label-cor]{Corresponding author.}

\affiliation[label-addr1]{
	organization={School of Information Science and Technology, ShanghaiTech University},
	city={Shanghai},
	postcode={201210},
	country={China.}}

\affiliation[label-addr2]{
	organization={College of Mathematics and Statistics, Chongqing University},
	city={Chongqing},
	postcode={401331},
	country={China.}}

\affiliation[label-addr3]{
	organization={Key Laboratory of Nonlinear Analysis and its Applications (Chongqing University), Ministry of Education},
	city={Chongqing},
	postcode={401331},
	country={China.}}

\begin{abstract}
	In Koopman operator theory, a finite-dimensional nonlinear system is transformed into an infinite but linear system using a set of observable functions. However, manually selecting observable functions that span the invariant subspace of the Koopman operator based on prior knowledge is inefficient and challenging, particularly when little or no information is available about the underlying systems. Furthermore, current methodologies tend to disregard the importance of the invertibility of observable functions, which leads to inaccurate results. To address these challenges, we propose the so-called FlowDMD, aka Flow-based Dynamic Mode Decomposition, that utilizes the Coupling Flow Invertible Neural Network (CF-INN) framework. FlowDMD leverages the intrinsically invertible characteristics of the CF-INN to learn the invariant subspaces of the Koopman operator and accurately reconstruct state variables. Numerical experiments demonstrate the superior performance of our algorithm compared to state-of-the-art methodologies.
\end{abstract}


\end{frontmatter}

\section{Introduction}
\label{sec-Introduction}

Nonlinear dynamic systems are widely prevalent in both theory and engineering applications. Since the governing equations are generally unknown in many situations, it can be challenging to study the systems directly based on the first principles. Fortunately, the data about the systems of interest could be available by experiments or observations. Instead, one could seek to understand the behavior of the nonlinear system through the data-driven approaches~\cite{brunton2016discovering, long2018pdenet, raissi2018deep, fuentes2021equation, kim2021integration}. 

The Koopman operator \cite{koopman1931hamiltonian}, which embeds the nonlinear system of interest into an infinite dimensional linear space by observable functions has attracted lots of attention. The Koopman operator acts on the infinite dimensional Hilbert space and aims to capture the full representations of the nonlinear systems. Dynamic mode decomposition (DMD) calculates the spectral decomposition of the Koopman operator numerically by extracting dynamic information from the collected data. Concretely, DMD devises a procedure to extract the spectral information directly from a data sequence without an explicit formulation of the Koopman operator, which is efficient for handling high dimensional data \cite{schmid2022dynamic}. Variants of DMD are proposed to address challenges in different scenarios \cite{Jonathan2014OnDynamic, jovanovic2014sparsity, takeishi2017bayesian, arbabi2017ergodic, le2017higher, erichson2019randomized, zhang2019online, colbrook2023residual}.

The selection of observable functions plays an essential role in the DMD algorithm. Exact DMD \cite{Jonathan2014OnDynamic} exploits the identity mapping as the observables. This implies that one uses a linear system to approximate a nonlinear system with given data \cite{kutz2016dynamic}. This would yield inaccurate or even completely mistaken outcomes. Furthermore, the short-term prediction of Exact DMD might be acceptable for some cases, but the long-term prediction is probably unreliable. Typically, prior knowledge is required to select the observable functions that span the invariant subspace of the Koopman operator. However, the invariant subspace is not simply available. In order to overcome the limitations of the Exact DMD algorithm and capture the full feature of the nonlinear system, several data-driven selection strategies for observable functions have been proposed. Extended DMD (EDMD) \cite{williams2015data} lifts the state variables from the original space into a higher dimensional space using the dictionary functions. The accuracy and rate of convergence of EDMD depend on the choice of the dictionary functions. Therefore, EDMD needs as many dictionary functions as possible. This implies that the set of dictionary functions (nonlinear transformations) should be sufficiently complex, which results in enormous computational cost. Kernel based DMD (KDMD) \cite{williams2015kernelbased} differs from EDMD in that it utilizes the kernel trick to exploit the implicit expression of dictionary functions, whereas EDMD uses the explicit expression of dictionary functions. Nonetheless, both EDMD and KDMD are prone to overfitting \cite{otto2019linearly}, which leads to large generalization error. How to efficiently choose the observable functions that span the invariant subspace of the Koopman operator becomes a significant challenge. 

In contrast to EDMD and KDMD, observable functions can be represented by neural networks. Dictionary learning~\cite{li2017extended} couples the EDMD with a set of trainable dictionary functions, where dictionary functions are represented by a fully connected neural network (FNN) and an untrainable component. Fixing the partial dictionary function facilitates the reconstruction of the state variables, however, this setting implicitly assumes that linear term lies in the invariant subspace of the Koopman operator. \citet{yeung2019learning} select low-dimensional dictionary functions more efficiently using deep neural networks. 

Autoencoder (AE) neural networks have been widely applied to learn the optimal observable functions and reconstruction functions in Koopman embedding \cite{otto2019linearly, takeishi2017learning, lusch2018deep, azencot2020forecasting, pan2020physicsinformed, LI2021110660}. Concretely, the invariant subspace of the Koopman operator and reconstruction functions are represented by the encoder and decoder network in AE, respectively. \citet{lusch2018deep} utilize neural networks to identify the Koopman eigenfunctions and introduced an auxiliary network to cope with the dynamic systems with continuous spectrum. \citet{azencot2020forecasting} propose the Consistent Koopman AE model that combines the forward-backward DMD method \cite{dawson2016characterizing} with the AE model. This approach extracts the latent representation of high-dimensional nonlinear data and eliminates the effect of noise in the data simultaneously. \citet{pan2020physicsinformed} parameterize the structure of the transition matrix in linear space and construct an AE model to learn the residual of the DMD. \citet{LI2021110660} utilize deep learning and the Koopman operator to model the nonlinear multiscale dynamical problems, where coarse-scale data is used to learn the fine-scale information through a set of multiscale basis functions. \citet{wang2023koopman} propose Koopman Neural Forecaster combining AE with Koopman operator theory to predict the data with distributional shifts.

Representing Koopman embedding by dictionary learning or AE networks has several drawbacks. Firstly, the reconstruction in dictionary learning partially fixes the dictionary functions, which leads to a low level of interpretability of the model. Secondly, the encoder and decoder in an AE model are trained simultaneously, but neither of them is invertible, cf. \cite{alford2022deep} for more details. Moreover, due to the structural noninvertibility of the encoder and decoder, it typically requires a large amount of training data in order to obtain accurate representations, which makes the AE model prone to overfitting. \citet{alford2022deep} analyze the property of both the encoder and decoder in AE and proposed the deep learning dynamic mode decomposition.  \citet{bevanda2022learning} constructed a conjugate map between the nonlinear system and its Jacobian linearization, which is learned by a diffeomorphic neural network.

In this paper, we develop a novel architecture called FlowDMD, aka Flow-based Dynamic Mode Decomposition, to learn the Koopman embedding. Specifically, we apply the coupling flow invertible neural networks to learn the observable functions and reconstruction functions. The invertibility of the learned observable functions makes our method more flexible than dictionary learning and AE learning. Our contributions are three-folds:
\begin{enumerate}
    \item The state reconstruction is accomplished by the backward direction (inversion) of the CF-INN, which increases the interpretability of the neural network and alleviates the overfitting of AE. 
    \item Due to the structural invertibility of CF-INN, the loss function for the state reconstruction is simplified compared with AE, which makes the network training easier.
    \item The parameters to be optimized are reduced dramatically since the learned mappings and their inverse share the same parameters.
\end{enumerate}

This paper is organized as follows. In Section~\ref{sec-Preliminaries}, we briefly review the Koopman operator theory and DMD.  In Section~\ref{sec-Learning}, we present the structure of CF-INN and introduce how to learn the invariant subspace of the Koopman operator and the reconstruction functions. In Section~\ref{sec-Numerical}, several numerical experiments are performed to demonstrate the performance of our method, and we summarize our work in Section~\ref{sec-Conclusion}.

\section{Preliminaries}
\label{sec-Preliminaries}

\subsection{Koopman operator theory}

Consider the nonlinear autonomous system in discrete form,
\begin{equation}\label{equ:nonlinear}
	\mathbf{x}_{k+1} = f(\mathbf{x}_k), \quad \mathbf{x}_k \in \mathcal{M} \subset \mathbb{R}^m,
\end{equation}
where $\mathcal{M}$ represents the set of state space, $f: \mathcal{M} \rightarrow \mathcal{M}$ is an unknown nonlinear map, and $k$ is the time index. 

\begin{definition}[Koopman operator \cite{kutz2016dynamic}]
	For the nonlinear system \eqref{equ:nonlinear}, the Koopman operator $\mathcal{K}$ is an infinite-dimensional linear operator that acts on all observable functions $g:\mathcal{M} \rightarrow \mathbb{C}$ such that
	\[
		\mathcal{K}g(\mathbf{x}) = g(f(\mathbf{x})).
	\]
Here, $g(x)\in\mathcal{H}$ and $\mathcal{H}$ represents the infinite dimensional Hilbert space.
\end{definition}

Through the observable functions, the nonlinear system~\eqref{equ:nonlinear} could be transformed into an infinite-dimensional linear system using the Koopman operator,
\begin{equation}\label{equ:linear}
	g(\mathbf{x}_{k+1}) = g(f(\mathbf{x}_k)) = \mathcal{K}g(\mathbf{x}_k).
\end{equation}
Note that the Koopman operator is linear, \textit{i.e.}, $\mathcal{K}(\alpha_1 g_1(\mathbf{x}) + \alpha_2 g_2(\mathbf{x})) = \alpha_1g_1(f(\mathbf{x})) + \alpha_2g_2(f(\mathbf{x}))$, with $g_1(\mathbf{x}), g_2(\mathbf{x})\in \mathcal{H}$ and $\alpha_1, \alpha_2 \in \mathbb{R}$. As $\mathcal{K}$ is an infinite-dimensional operator, we denote its eigenfunctions and eigenvalues by $\{\lambda_i, \varphi _{i}( x)\}_{i=0}^{\infty }$ such that $\mathcal{K} \varphi _{i}(\mathbf{x}) =\lambda _{i} \varphi _{i}(\mathbf{x})$, where $\varphi _{i}(\mathbf{x}) :\mathcal{M}\rightarrow \mathbb{R}$, $\lambda _{i} \in \mathbb{C}$. 

The Koopman eigenfunctions define a set of intrinsic measurement coordinates, then a vector-valued observable function $\mathbf{g}(\mathbf{x}) = [g_1(\mathbf{x}), \cdots, g_n(\mathbf{x})]^T$ could be written in terms of the Koopman eigenfunctions, 
\begin{equation}\label{equ:decom-1}
	\mathbf{g} (\mathbf{x}_{k} )=\begin{bmatrix}
		g_{1}(\mathbf{x}_{k})\\
		\vdots \\
		g_{n}(\mathbf{x}_{k})
		\end{bmatrix} =\sum _{i=1}^{\infty } \varphi _{i}(\mathbf{x}_{k})\begin{bmatrix}
		< \varphi _{i} ,g_{1}> \\
		\vdots \\
		< \varphi _{i} ,g_{n}> 
	\end{bmatrix} =\sum _{i=1}^{\infty } \varphi _{i} (\mathbf{x}_{k} )\mathbf{v}_{i},
\end{equation}
where $\mathbf{v}_{i}$ refers to the $i$-th Koopman mode with respect to the Koopman eigenfunction $\varphi _{i}(\mathbf{x})$. Combining \eqref{equ:linear} and \eqref{equ:decom-1}, we have the decomposition of a vector-valued observable functions
\begin{equation*}\label{equ:decom-2}
	\mathbf{g}(\mathbf{x}_{k+1}) =\mathcal{K} \mathbf{g}(\mathbf{x}_{k}) =\mathcal{K}\sum _{i=1}^{\infty} \varphi _{i} (\mathbf{x}_{k} )\mathbf{v}_{i} =\sum _{i=1}^{\infty} \lambda _{i} \varphi _{i} (\mathbf{x}_{k} )\mathbf{v}_{i}.
\end{equation*}
Furthermore, the decomposition could be rewritten as 
\begin{equation*}\label{equ:decom-3}
	\mathbf{g}(\mathbf{x}_{k}) =\sum _{i=1}^{\infty} \lambda _{i}^{k} \varphi _{i} (\mathbf{x}_{0} )\mathbf{v}_{i}.
\end{equation*}

In practice, we need a finite-dimensional representation of the infinite-dimensional Koopman operator. Denote the $n$-dimensional invariant subspace of the Koopman operator $\mathcal{K}$ by $\mathcal{H}_{g}$, \textit{i.e.}, $\forall g(\mathbf{x}) \in \mathcal{H}_{g}, \mathcal{K}g(\mathbf{x}) \in \mathcal{H}_{g}$. Let $\{g_{i}(\mathbf{x})\}_{i=1}^{n}$ be one set of basis of $\mathcal{H}_{g}$, this induces a finite-dimensional linear operator $\mathbf{K}$ \cite{kutz2016dynamic}, which projects the Koopman operator $\mathcal{K}$ onto $\mathcal{H}_{g}$, \textit{i.e.}, for the $n$-dimensional vector-valued observable functions $\mathbf{g}(\mathbf{x}) = [g_1(\mathbf{x}), \cdots, g_n(\mathbf{x})]^T$, we have 
\begin{equation}\label{equ:linear-finite}
	\mathbf{g} (x_{k+1} )=\begin{bmatrix}
		g_{1}( x_{k+1})\\
		\vdots \\
		g_{n}( x_{k+1})
		\end{bmatrix} =\begin{bmatrix}
		\mathcal{K}g_{1}( x_{k})\\
		\vdots \\
		\mathcal{K}g_{n}( x_{k})
		\end{bmatrix} =\mathbf{K}\begin{bmatrix}
		g_{1}( x_{k})\\
		\vdots \\
		g_{n}( x_{k})
	\end{bmatrix} =\mathbf{K}\mathbf{g} (x_{k} )
\end{equation}

\subsection{Dynamic mode decomposition}
DMD approximates the spectral decomposition of the Koopman operator numerically. Given the state variables $\displaystyle \{\mathbf{x}_0, \mathbf{x}_1,\cdots, \mathbf{x}_p\}$ and a vector-valued observable function $\mathbf{g}(\mathbf{x}) = [g_1(\mathbf{x}), \cdots, g_n(\mathbf{x})]^T$, then we get the sequence $\displaystyle \{\mathbf{g}(\mathbf{x}_0), \mathbf{g}(\mathbf{x}_1),\cdots, \mathbf{g}(\mathbf{x}_p)\}$, where each $\mathbf{g}(\mathbf{x}_k)\in\mathbb{R}^n$ is the observable snapshot of the $k$-th time step. According to \eqref{equ:linear-finite}, we have 
\begin{equation*}\label{back-equ-1}
	\mathbf{g}(\mathbf{x}_{k+1}) = \mathbf{K} \mathbf{g}(\mathbf{x}_{k}),
\end{equation*}
where $\mathbf{K}\in\mathbb{R}^{n\times n}$ is the matrix form of the finite-dimensional operator. 
For the two data matrices, $\mathbf{X}_\mathbf{g}=[\mathbf{g}(\mathbf{x}_0), \cdots, \mathbf{g}(\mathbf{x}_{p-1})]$ and $\mathbf{Y}_\mathbf{g}=[\mathbf{g}(\mathbf{x}_1), \cdots, \mathbf{g}(\mathbf{x}_{p})]$, where $\mathbf{X}_\mathbf{g}$ and $\mathbf{Y}_\mathbf{g}$ are both in $\mathbb{R}^{n\times p}$, which satisfies $\mathbf{Y}_\mathbf{g}=\mathbf{K}\mathbf{X}_\mathbf{g}$. Therefore, $\mathbf{K}$ can be represented by
\begin{equation*}
	\mathbf{K}=\mathbf{Y}_\mathbf{g}\mathbf{X}^\dagger_\mathbf{g},
\end{equation*}
where $\mathbf{X}^\dagger_\mathbf{g}$ denotes the Moore-Penrose inverse of $\mathbf{X}_\mathbf{g}$. 

The Exact DMD algorithm developed by \citet{Jonathan2014OnDynamic} computes dominant eigen-pairs (eigenvalue and eigenvector) of $\mathbf{K}$ without the explicit formulation of $\mathbf{K}$. In Algorithm \ref{alg:ExactDMD}, we present the DMD algorithm on the observable space, which is a general form of the Exact DMD algorithm. When using the identical mapping as the observable functions, \textit{i.e.}, $\mathbf{g}(\mathbf{x}) = \mathbf{x}$,  Algorithm \ref{alg:ExactDMD} is identical to the Exact DMD algorithm.

\begin{algorithm}[htbp]
	\caption{DMD on observable space \cite{kutz2016dynamic, lu2020prediction}}
	\label{alg:ExactDMD}
	\begin{algorithmic}
	\STATE{1. Compute the (reduced) SVD of $\mathbf{X}_\mathbf{g}$, $\mathbf{X}_\mathbf{g}=\mathbf{U_r}\mathbf{\Sigma_r}\mathbf{V_r^*}$, where $\mathbf{U_r}\in\mathbb{C}^{n \times r}$, $\mathbf{\Sigma_r} \in\mathbb{R}^{r \times r}$, $\mathbf{V_r}\in\mathbb{C}^{p \times r}$.}
	\STATE{2. Compute $\tilde{\mathbf{K}}=\mathbf{U_r^*}\mathbf{Y}_\mathbf{g}\mathbf{V_r}\mathbf{\Sigma_r^{-1}}$.}
	\STATE{3. Compute the eigen-pairs of $\tilde{\mathbf{K}}$: $\tilde{\mathbf{K}}\mathbf{W}=\mathbf{W\Lambda}$.}
	\STATE{4. Reconstruct the eigen-pairs of $\mathbf{K}$, where eigenvalues of $\mathbf{K}$ are diagonal entries of $\Lambda$, the corresponding eigenvectors of $\mathbf{K}$(DMD modes) are columns of $\mathbf{\Phi} = \mathbf{Y}_\mathbf{g}\mathbf{V_r}\mathbf{\Sigma_r^{-1}}\mathbf{W}$.}
	\STATE{5. Approximate the observation data via DMD, $\hat{\mathbf{g}}(\mathbf{x}_{k}) = \mathbf{\Phi} \mathbf{\Lambda}^{k}\mathbf{b}$, where $\mathbf{b}=\mathbf{\Phi}^{\dagger}\mathbf{g}(\mathbf{x}_0)$.}
	\STATE{6. Reconstruct the state variables $\hat{\mathbf{x}}_{k} =\mathbf{g}^{-1}(\hat{\mathbf{g}}(\mathbf{x}_{k})) =\mathbf{g}^{-1}\left(\mathbf{\Phi \Lambda }^{k}\mathbf{b}\right)$.}
	\end{algorithmic}
\end{algorithm}

\subsection{State reconstruction}

Koopman operator theory utilizes observable functions $\mathbf{g}$ to transform the nonlinear system \eqref{equ:nonlinear} into a linear system while preserving the nonlinearity. Evolving the nonlinear system \eqref{equ:nonlinear} is computationally expensive or even impossible when $f$ is unknown, whereas evolving through the Koopman operator \eqref{equ:linear} offers a promising and computationally efficient approach.

\begin{figure}[ht]
	\centering
	\includegraphics[width=0.8\textwidth]{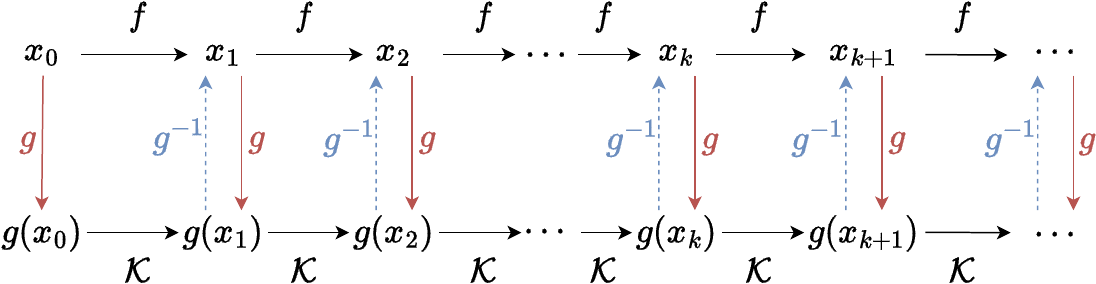}
	\caption{Koopman operator and inverse of observable functions}
	\label{fig:koopman-fig1}
\end{figure}

Figure \ref{fig:koopman-fig1} illustrates the relation between the nonlinear evolution $f$ and the Koopman operator evolution where the system evolves linearly in the observation space $\mathcal{H}$. By computing the Koopman eigenvalues and modes, we can make predictions of the observable functions $\mathbf{g}(\mathbf{x})$. We could reconstruct the state $\mathbf{x}$ by the inverse of the observable functions $\mathbf{g}^{-1}(\mathbf{x})$ provided that $\mathbf{g}(\mathbf{x})$ is invertible. The invertibility of observable functions is essential to ensure the reconstruction accuracy and the interpretability of the outcomes.

\vspace{-0.2cm}
Typical observable functions $\mathbf{g}(\mathbf{x})$ selection are performed manually based on prior knowledge. 
Exact DMD takes the identical mapping, 
while the EDMD utilizes a set of pre-defined functions such as polynomials, Fourier modes, radial basis functions, and so forth \cite{williams2015data}. However, these methods can be inaccurate and inefficient for Koopman embeddings learning. Deep neural networks, as efficient global nonlinear approximators, could be applied to represent the observable function $\mathbf{g}(\mathbf{x})$ and the reconstruction function $\mathbf{g}^{-1}(\mathbf{x})$.
Several studies have demonstrated that the encoder and decoder networks in AE correspond to $\mathbf{g}(\mathbf{x})$ and $\mathbf{g}^{-1}(\mathbf{x})$, respectively \cite{otto2019linearly, takeishi2017learning,  lusch2018deep, azencot2020forecasting, pan2020physicsinformed, LI2021110660}.

In practical applications, it is not always guaranteed that $\mathbf{g}(\mathbf{x})$ is invertible. In the learning Koopman embedding via AE, the invertibility of $\mathbf{g}(\mathbf{x})$ is enforced through numerical constraints, \textit{i.e.}, the reconstruction error $\|\mathbf{x} - \mathbf{g}^{-1}(\mathbf{g}(\mathbf{x}))\|_2^2$, which tends to result in overfitting and suboptimal performance \cite{alford2022deep}. Besides, the reconstruction error is trained simultaneously with the prediction error and the linearity error \cite{lusch2018deep}. The weights assigned to each loss term are hyperparameters that can be challenging to tune. In this paper, we propose a structurally invertible mapping learning framework, which eliminates the need for the reconstruction term in the loss function and yields more robust and accurate results. We present the details of our method in Section \ref{sec-Learning}.

\section{Learning Koopman embedding by invertible neural networks}
\label{sec-Learning}
In this section, we first briefly review the AE neural network and demonstrate the limitation of this class of neural networks in the Koopman embedding learning. Then, we introduce our method to overcome this limitation. 

For notational simplicity, we introduce some notations herein. For two mappings or functions $f_1$ and $f_2$, their composite $f_2(f_1(x))$ is denoted by $f_2 \circ f_1(x)$. Given two vectors $\mathbf{a} = [a_1, \cdots, a_m]^T, \mathbf{b} = [b_1, \cdots, b_m]^T \in \mathbb{R}^{m}$, their Hadamard product is the element-wise multiplication, represented by $\mathbf{a} \odot \mathbf{b} = [a_1b_1, a_2b_2, \cdots, a_mb_m]^T$. Identically, their Hadamard division is defined as the element-wise division, \textit{i.e.}, $\mathbf{a} \oslash \mathbf{b} = [a_1/b_1, a_2/b_2, \cdots, a_m/b_m]^T$.
\subsection{Drawback of AE in the Koopman embedding learning}
Most of the work use the AE neural networks as the backbone to learn the invariant subspace of the Koopman operator and reconstruct the state variables. AE as the frequently-used unsupervised learning structure of neural networks, consists of two parts, \textit{i.e.}, the encoder $\mathcal{E}$ and the decoder $\mathcal{D}$. 
AE learns these two mappings (functions) $\mathcal{E}$ and $\mathcal{D}$ by optimizing

\[
\min_{\mathcal{E, D}} \mathbb{E}_{x\sim m( x)}[\mathrm{loss}( x, \mathcal{D\circ E}( x))].
\]
Here $m(x)$ denotes the distribution of the input data, $\text{loss}(x,y)$ describes the difference between $x$ and $y$, and $\mathbb{E}(\cdot)$ represents the expectation.

\begin{definition}\label{def:inv}
    Let $f_1:S \rightarrow S'$ be an arbitrary mapping, and it is said to be invertible if there exists a mapping $f_2: S' \rightarrow S$ such that
    \[
    f_1\circ f_2 = \mathcal{I}, f_2\circ f_1 = \mathcal{I},
    \]
    where $\mathcal{I}$ is the identity mapping. Then,  $f_2$ is said to be the inverse mapping of $f_1$.
\end{definition}

Let $\mathcal{E}$ and $\mathcal{D}$ be two mappings learned by AE such that $\mathcal{D} \circ \mathcal{E}\approx\mathcal{I}$. However, the reverse order of the mapping $\mathcal{E} \circ \mathcal{D}$ is not always a good approximation to the identity mapping, moreover, $\mathcal{E}$ and $\mathcal{D}$ are generally not invertible \cite{alford2022deep}. The main reason is that while AE strives to reach $\mathcal{D}\circ \mathcal{E}\approx \mathcal{I}$, it omits the additional constraint $\mathcal{E}\circ\mathcal{D}\approx\mathcal{I}$ which requires the latent variable data to train. Unfortunately, the latent variables are not accessible, thus rendering it impossible for AE to satisfy $\mathcal{E}\circ\mathcal{D}\approx\mathcal{I}$ and $\mathcal{D}\circ\mathcal{E}\approx\mathcal{I}$ simultaneously.

\begin{figure}[!t]
    \begin{center}
        \subfloat[]{
            \includegraphics[width=0.45\linewidth]{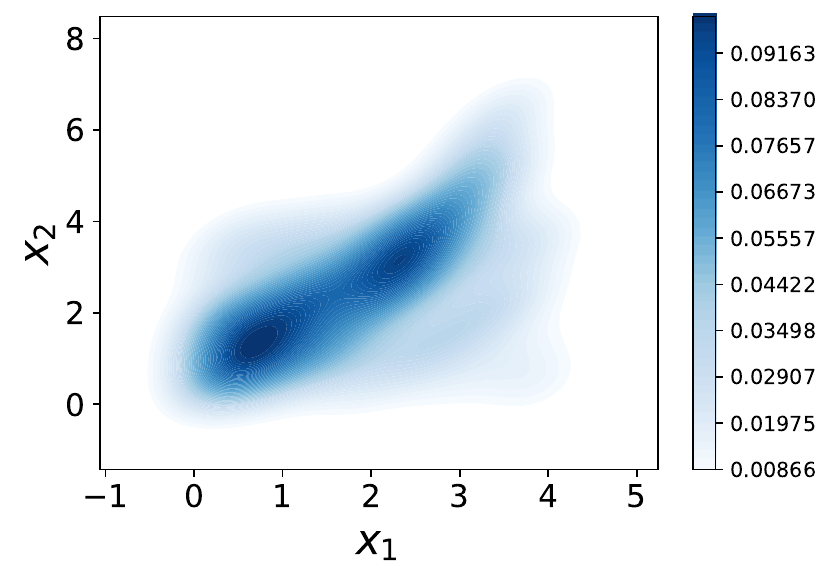}
        }
        \subfloat[]{
            \includegraphics[width=0.45\linewidth]{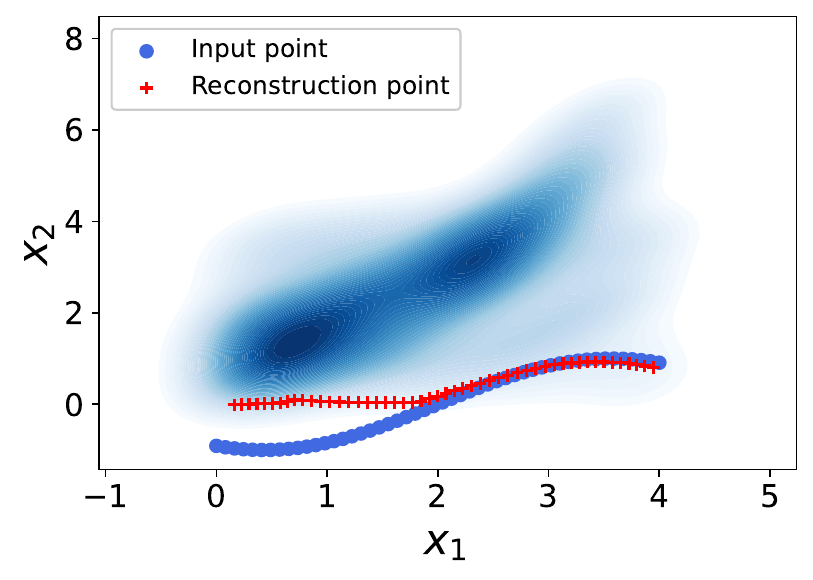}
        }\\
        \vspace{-.4cm}
        \subfloat[]{
            \includegraphics[width=0.45\linewidth]{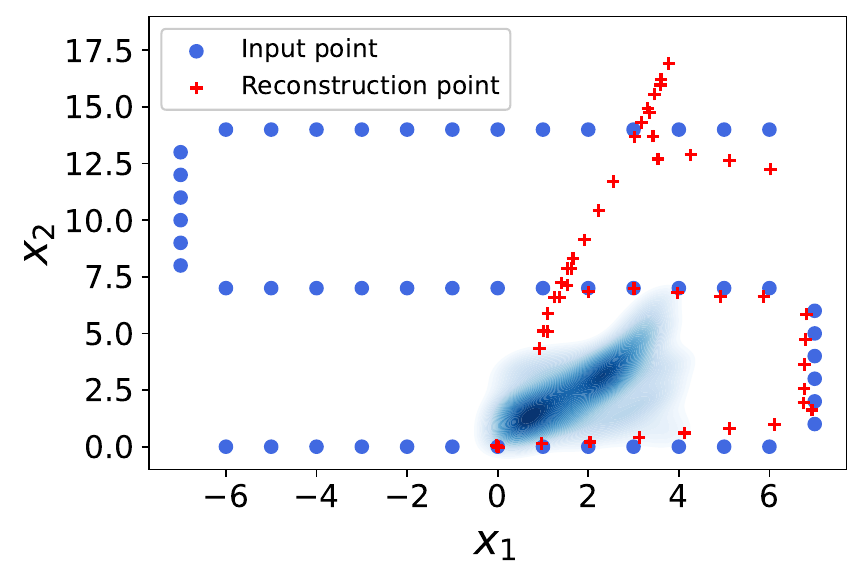}
        }
        \subfloat[]{
            \includegraphics[width=0.45\linewidth]{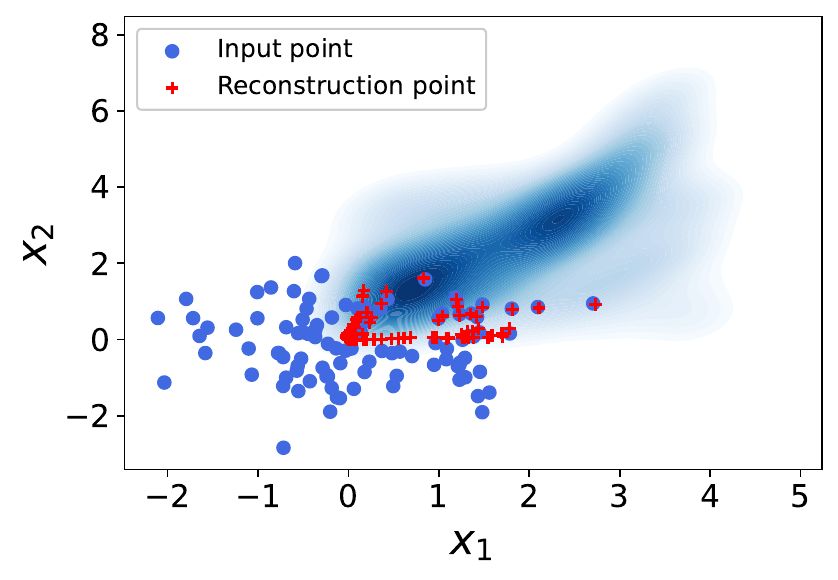}
        }
    \end{center}
    \vspace{-.4cm}
    \caption{Generalization capability test of AE. 
    (a) the training data distribution.
    (b)  the $sin(x)$ test function.
    (c) S-shaped scatters test.
    (d) random scatters from 2-d standard normal distribution.}\label{fig:comp1}
\end{figure}

AE learns an identity mapping $\mathcal{I}$ from a training data set $\mathcal{S}$, \textit{i.e.}, for any $\mathbf{x}\in \mathcal{S}, \mathcal{D}\circ \mathcal{E}(\mathbf{x})\approx \mathbf{x}$. For data out of the set $\mathcal{S}$, the mapping learned by AE may perform badly. In other words, AE may have poor generalization capability. Next, we use a preliminary experiment to demonstrate this limitation. The details of this numerical example are given in Section \ref{sec:fixed_point}.
We use the structure of AE defined in \cite{LI2021110660} and randomly generate 120 trajectories to train the AE, and the results are depicted by Figure \ref{fig:comp1}.

Figure \ref{fig:comp1} compares the input data points out of the distribution of the training data with the corresponding reconstructed data points using  the trained AE model. Figure \ref{fig:comp1}(a) shows the density distribution of training data set $\mathcal{S}$, which provides a rough illustration of the data space $\mathcal{S}$.  For  the reconstruction test of AE, we generate three types of data, \textit{i.e.}, the sin-shaped scatters, the S-shaped scatters, and scatters from the standard 2-d normal distribution.  We plot the corresponding input points (blue) and reconstructed data points (red) of the AE. The results shown in the next three subfigures illustrate that  AE can reconstruct the input data points nearby the  training data set $\mathcal{S}$ very well. But for the data points far away from $\mathcal{S}$, AE performs badly. The same situation happens in learning the Koopman embedding. Specifically, in the training process of AE, one aims to find the Koopman invariant space by minimizing the error of the  Koopman embedding learning and the reconstruction error. However, minimizing the error between latent variables and their corresponding reconstruction denoted by $\text{loss}(\mathbf{x},\mathcal{E}\circ\mathcal{D}(\mathbf{x}))$ is intractable. This result is in poor stability and generalization capability.

\subsection{Structure of CF-INN}
We have shown that the mapping learned by AE performs poorly, which inspires us that invertibility can greatly reduce computational complexity and yields better generalization capability.
Next, we introduce an invertible neural network to overcome the drawback of AE.
Let $\mathbf{g}_{\theta}(\mathbf{x}):\mathbf{X}\rightarrow \mathbf{Y}$ denote the input-output mapping of the invertible neural network, where ${\theta}$ represents the parameters of the neural network. Let $\mathbf{f}_{\theta}$ be the inverse mapping of $\mathbf{g}_{{\theta}}$ which shares the same parameters with $\mathbf{g}_{{\theta}}$. Then we can reconstruct $\mathbf{x}$ in the backward direction by $ \mathbf{f}_{\theta}(\mathbf{y}):\mathbf{Y} \rightarrow \mathbf{X}$. 
In generative tasks of machine learning, the forward generating direction is called the flow direction and the backward direction is called the normalizing direction. Next, we introduce the concept of coupling flow (CF), which belongs to the category of invertible neural networks.

\begin{definition}[Coupling flow~\cite{papamakarios2021normalizing}]\label{def:CF} 
  Let $m \in \mathbb{N}$ and $m \geq 2$, we partition a vector $\mathbf{z} \in \mathbb{R}^m$ as $\mathbf{z} = \begin{bmatrix}
    \mathbf{z}_{up}\\
    \mathbf{z}_{low}
    \end{bmatrix}$ with $\mathbf{z}_{up} \in \mathbb{R}^q$ and $\mathbf{z}_{low} \in \mathbb{R}^{m-q}$ for $1\leq q \leq m-1$.
    The coupling flow $h_{q, \tau}: \mathbb{R}^{m} \rightarrow \mathbb{R}^{m}$ is defined by
    \begin{equation}\label{equ:fun-h}
        h_{q, \tau}(\mathbf{z}) =\begin{bmatrix}
            \mathbf{z}_{up}\\
            \tau (\mathbf{z}_{low}, \sigma(\mathbf{z}_{up}))
            \end{bmatrix},
    \end{equation}
    where $\sigma: \mathbb{R}^q \rightarrow \mathbb{R}^l$ is an arbitrary mapping, and $\tau(\cdot, \sigma(\mathbf{y})): \mathbb{R}^{m-q} \times \mathbb{R}^l \rightarrow \mathbb{R}^{m-q}$ is a bijective mapping for any $\mathbf{y}\in\mathbb{R}^q$.
\end{definition}

The CF given by Definition~\ref{def:CF} is invertible (bijective) if and only if $\tau$ is bijective~\cite{9089305}. 
Specifically, let $\mathbf{w} = h_{q, \tau}(\mathbf{z})$ and partition it in the same manner with $\mathbf{z}$, \textit{i.e.}, $\mathbf{w} = \begin{bmatrix}
    \mathbf{w}_{up}\\
    \mathbf{w}_{low}
    \end{bmatrix}$, where $\mathbf{w}_{up} \in \mathbb{R}^{q}$, and $\mathbf{w}_{low} \in \mathbb{R}^{m-q}$. 
Then we can obtain the inverse of~\eqref{equ:fun-h} given by, 
    \[
    \mathbf{z} := h_{q, \tau}^{-1}(\mathbf{w})=
    \begin{bmatrix}
    \mathbf{w}_{up}\\
    \tau^{-1} (\mathbf{w}_{low}, \sigma(\mathbf{w}_{up}))
    \end{bmatrix}.
    \]
One of the mostly used  CF is the affine coupling flow (ACF)~\cite{dinh2014nice,dinh2016density,kingma2018glow}, where $\tau$ is an element-wise invertible (bijective) function. 
\begin{definition}[Affine coupling flow \cite{9089305}]\label{def:ACF}
Given $\mathbf{z} = \begin{bmatrix}
        \mathbf{z}_{up}\\
        \mathbf{z}_{low}
        \end{bmatrix} \in \mathbb{R}^m$ with $\mathbf{z}_{up} \in \mathbb{R}^q$ and $\mathbf{z}_{low} \in \mathbb{R}^{m-q}$ for $1\leq q \leq m-1$, the affine coupling flow $\psi_{q, s, t}: \mathbb{R}^m \rightarrow \mathbb{R}^{m}$ is defined by
    \begin{equation}\label{equ:ACF-fun}
        \psi_{q, s, t}(\mathbf{z}) =\begin{bmatrix}
        \mathbf{z}_{up}\\
        (\mathbf{z}_{low}+t(\mathbf{z}_{up})) \odot s(\mathbf{z}_{up})
        \end{bmatrix},
    \end{equation}
    where $s, t: \mathbb{R}^q \rightarrow \mathbb{R}^{m-q}$ are two arbitrary mappings.
\end{definition}

Equation~\eqref{equ:ACF-fun} is often referred to as the forward direction computations of ACF. Let $\mathbf{w} = \psi_{q, s, t}(\mathbf{z})$, we can give its corresponding backward direction computations by 
\[
\mathbf{z} := \psi_{q, s, t}^{-1}(\mathbf{w})=\begin{bmatrix}
    \mathbf{w}_{up}\\
    \mathbf{w}_{low}\oslash s(\mathbf{w}_{up})-t(\mathbf{w}_{up})
    \end{bmatrix},
    \] 
where $\mathbf{w}$ is partitioned in the same manner with $\mathbf{z}$. Additionally, the mappings $s$ and $t$ in Definition~\ref{def:ACF} can be any nonlinear functions or neural networks such as FNN.

Let $\psi_1,\dots,\psi_L$ be a sequence of ACFs and define $\mathbf{g}_{\theta}=\psi_L\circ \psi_{L-1} \circ \cdots \psi_{1}$, where ${\theta}$ represents the parameters of $\{\psi_i\}_{i=1}^L$. Thus, $\mathbf{g}_{\theta}$ results in an invertible neural network and is called by CF-INN in this paper. Moreover, the division index $q$ of the input vector $\mathbf{z}$ is user-guided. In this paper, we set $q=\lceil m/2 \rceil$, where $\lceil \cdot \rceil$ is the rounding function. Furthermore, in order to mix the information propagated in the network sufficiently, we could define a flipped ACF denoted by $\displaystyle\widetilde{\psi}_{q, s, t}: \mathbb{R}^m \rightarrow \mathbb{R}^{m}$ which is obtained by simply flipping two input parts of ACF. The forward direction computation of a flipped ACF is given by
\begin{equation*}
    \widetilde{\psi}_{q, s, t}(\mathbf{z})=\begin{bmatrix}
        (\mathbf{z}_{up}+t(\mathbf{z}_{low})) \odot s(\mathbf{z}_{low})\\
        \mathbf{z}_{low}
        \end{bmatrix}.
\end{equation*}

\begin{figure}[htbp]
    \begin{center}
        \includegraphics[width=.7\linewidth,trim=0 0 0 0]{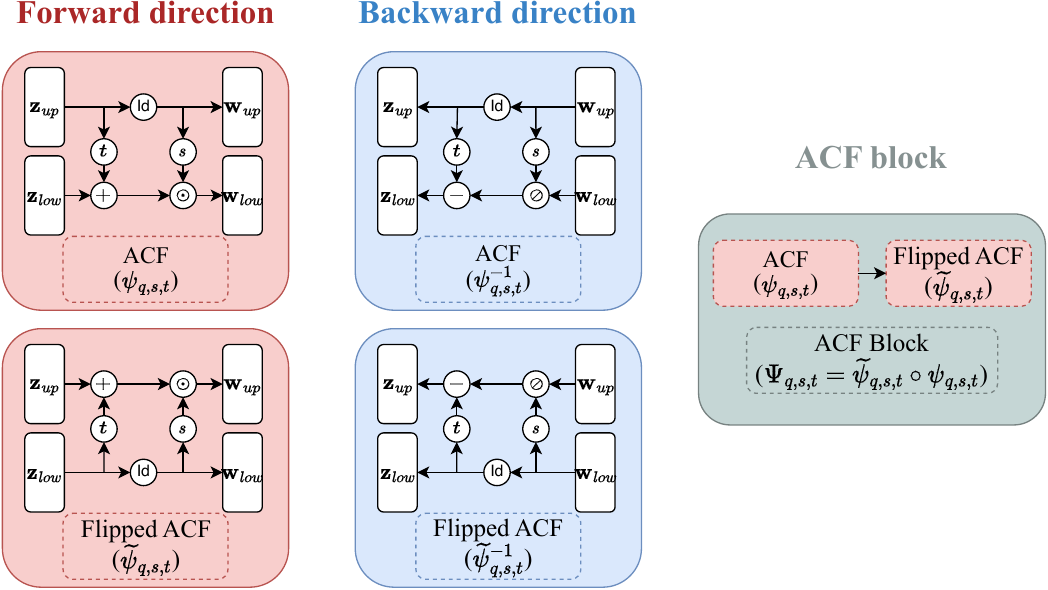}
    \end{center}
    \vspace{-.3cm}
    \caption{The forward and backward directions of ACF and flipped ACF, as well as the structure of an ACF block. Here, the ``Id'' operation represents the identity mapping.}\label{fig:Method}
\end{figure}

We can compose an ACF block denoted by $\Psi_{q, s, t} = \widetilde{\psi}_{q, s, t} \circ \psi_{q, s, t}$ using a standard ACF and a flipped ACF. The structure of an ACF $\psi_{q, s, t}$, a flipped ACF $\widetilde{\psi}_{q, s, t}$, and an ACF block $\Psi_{q, s, t}$ are shown in Figure~\ref{fig:Method} , where the left and middle columns represent the forward and backward computations of an ACF and a flipped ACF, respectively. The right column shows the structure of an ACF block, which is a CF-INN of depth 2.

When the depth of a CF-INN model, \textit{i.e.}, $L$, is large, its training becomes challenging. The main curse is that the dividend term  $s$ is too small in $\Psi_{q, s, t}$ in the backward direction computations. This issue can be solved by replacing the ACF with a residual coupling flow (RCF). Similar idea has also been applied in the residual term of ResNet.
\begin{definition}[Residual coupling flow~\cite{gomez2017reversible}]\label{def:RCF}
    Given $\mathbf{z} = \begin{bmatrix}
        \mathbf{z}_{up}\\
        \mathbf{z}_{low}
        \end{bmatrix} \in \mathbb{R}^m$ with $\mathbf{z}_{up} \in \mathbb{R}^q$ and $\mathbf{z}_{low} \in \mathbb{R}^{m-q}$ for $1\leq q \leq m-1$, the residual coupling flow $\psi_{q, t}: \mathbb{R}^m \rightarrow \mathbb{R}^m$ is defined by
    \[
    \psi_{q, t}(\mathbf{z}) =\begin{bmatrix}
    \mathbf{z}_{up}\\
    \mathbf{z}_{low}+t(\mathbf{z}_{up})
    \end{bmatrix},
    \]
    where $t: \mathbb{R}^q \rightarrow \mathbb{R}^{m-q}$ is an arbitrary mapping.
\end{definition}
RCFs are simplifications of ACFs and similar with an ACF block, we can obtain a RCF block in composition of a RCF and a flipped RCF, which is a simplified ACF block.

\subsection{Loss function of FlowDMD for Koopman embedding}
In this paper, we use the CF-INN to learn the Koopman invariant subspace and the reconstructions simultaneously, where the forward direction of CF-INN is represented by $\mathbf{g}_{\theta}$ and its backward direction is represented by $\mathbf{f}_{\theta}$. Our method is called FlowDMD as it integrates CF-INN and DMD to compute the finite dimensional Koopman operator approximation and reconstruct system states.

It is noteworthy that the dimensions of input and output of CF-INN are inherently the same, which implies that the system states and the Koopman invariant subspace share the same dimension for FlowDMD. This does not generally hold true for the Koopman operator learning. However, this does not simply imply that FlowDMD fails for the general cases, vice versa, our approach is generally applicable. Next, we discuss why this holds true.
    
\begin{figure}[ht]
    \centering
    \includegraphics[width=\linewidth]{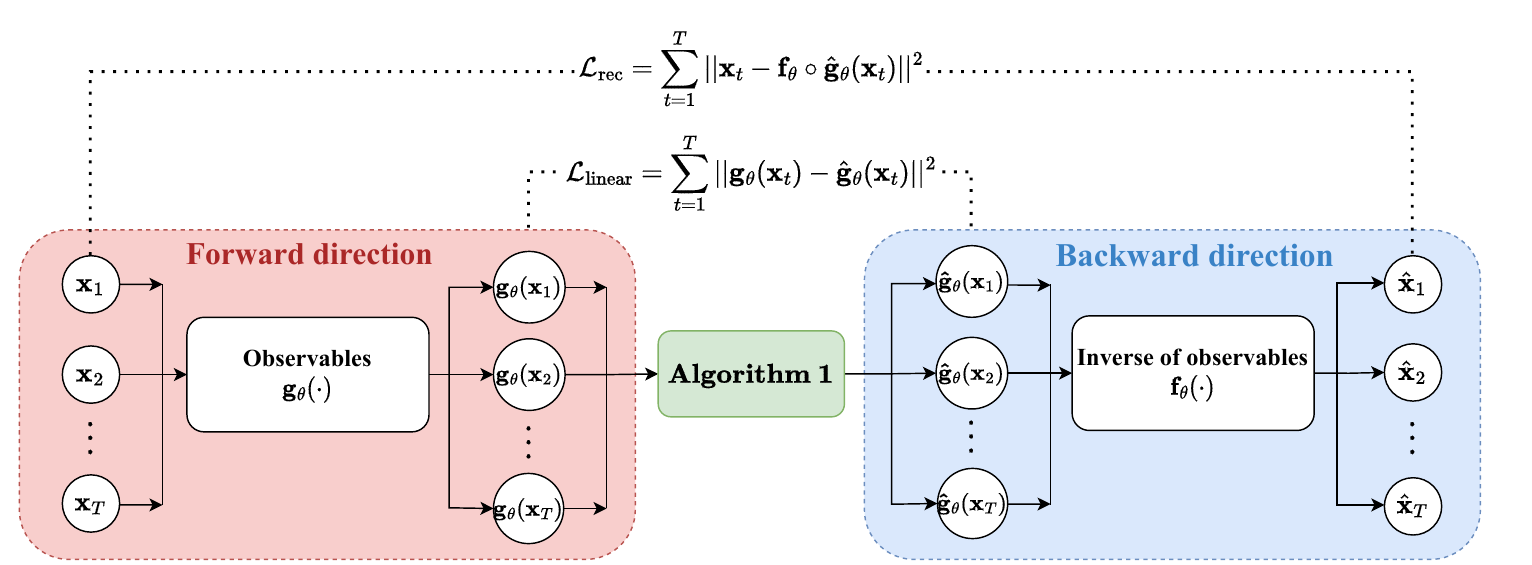}
    \caption{The general framework of FlowDMD. }\label{fig:total_framework}
\end{figure}

Consider the general case that dimension of the states being $m$ and the dimension of the Koopman invariant subspace being $\tilde{m}$.
\begin{enumerate}
    \item Case 1: $m\geq \tilde{m}$, the output of CF-INN is an $m$-dimensional vector functions which in turn gives an $m$-dimensional Koopman invariant subspace that already contains the $\tilde{m}$-dimensional Koopman invariant subspace. One can directly perform computations in this $m$-dimensional subspace which gives more accurate results than computing in the $\tilde{m}$-dimensional subspace without any extra cost. If strictly restricted to the $\tilde{m}$-dimensional Koopman invariant subspace, one can first project from the $m$-dimensional subspace to this subspace then perform computations and finally project back to the $m$-dimensional subspace for the state reconstruction using the backward direction of CF-INN. This procedure applies at the beginning and end of Algorithm 1 in Figure~\ref{fig:total_framework} and no modifications of CF-INN is needed. Note that the Koopman operator theory transforms nonlinear systems to linear systems by lifting system dimension which usually gives $m < \tilde{m}$. \\
    \item Case 2: $m < \tilde{m}$, we can augment the states by appending at least $\tilde{m}-m$ zeros in total, either (both) preceding or (and) succeeding the original states, which can be represented by $\mathbf{x}' = [0, \cdots, 0, \mathbf{x}^T, 0, \cdots, 0]^T\in\mathbb{R}^{m'}$. Using this simple technique, the CF-INN can be directly applied without any adjustment or even modifications of the loss function. The reconstructed states $\mathbf{\tilde{x}}' = [0, \cdots, 0, \mathbf{\tilde{x}}^T, 0, \cdots, 0]^T\in\mathbb{R}^{m'}$ has the same pattern with $\mathbf{x}'$ by prescription. Such methodology is analogous to the zero-padding technique commonly employed by image processing.
\end{enumerate}

The loss function of FlowDMD has two components which consists of the DMD approximation error and the state reconstruction error. Firstly, the observable functions evolve linearly in the Koopman invariant subspace. Hence, the linearity constrained loss function that represents the DMD approximation error is given by
\begin{equation*}
    \mathcal{L}_{\text{linear}}=\sum_{t=1}^T || \mathbf{g}_{\theta}(\mathbf{x}_t)-\hat{\mathbf{g}}_{\theta}(\mathbf{x}_t)||^2,
\end{equation*}
where $\hat{\mathbf{g}}_{\theta}(\mathbf{x}_t) = \Phi\Lambda^t \Phi^{\dagger}\mathbf{g}_{\theta}(\mathbf{x}_0)$ is the DMD approximation of the observable functions $\{\mathbf{g}(\mathbf{x}_t)\}_{t=1}^T$ by using Algorithm~\ref{alg:ExactDMD}. 

Secondly, the inverse mapping of $\mathbf{g}$, \textit{i.e}, the backward direction of CF-INN, $\mathbf{f}_\theta$, are used to reconstruct the states $\mathbf{x}_t$. Here, $\mathbf{f}_\theta$ shares the same network structure and parameters with $\mathbf{g}_\theta$. Therefore, the computational cost is greatly reduced, compared with AE that another neural network is required to parameterize the inverse mapping of $\mathbf{g}_\theta$. The reconstruction loss due to the DMD approximation error is given by
\begin{equation*}
    \mathcal{L}_{\text{rec}} =\sum _{t=1}^{T} ||\mathbf{x}_{t} -\mathbf{f}_{\theta} \circ \hat{\mathbf{g}}_{\theta} (\mathbf{x}_{t} )||^{2}.
\end{equation*}
The optimal parameters ${\theta}^*$ is determined by 
\begin{align*}
    {\theta}^*&= \mathop{\arg\min}\limits_{{\theta}} \mathcal{L}_{\text{linear}} + \alpha \mathcal{L}_{\text{rec}}, 
\end{align*}
where $\alpha$ is a user-guard hyperparameter. 

We summarize our FlowDMD framework for the Koopman embedding learning in Figure \ref{fig:total_framework}. In other Koopman embedding learning frameworks \cite{lusch2018deep, LI2021110660}, the reconstruction error induced by the noninvertibility of neural networks denoted by $\mathcal{L}_{\text{nn}} =\sum _{t=1}^{T} ||\mathbf{x}_{t} -\mathbf{f}_{\theta}(\mathbf{g}_{\theta} (\mathbf{x}_{t} ))||^{2}$ also needs to be considered. However, in our model, this term vanishes due to the invertibility of CF-INN, resulting in a notably simplified loss function, which makes the network training easier compared with~\cite{lusch2018deep, LI2021110660}.

\section{Numerical experiments}
\label{sec-Numerical}
In this section, we use three numerical examples to demonstrate the efficiency of our method for learning the Koopman embedding and compare its performance with LIR-DMD \cite{LI2021110660}, Exact DMD, and  EDMD. We use the Python library \textit{FEniCS}~\cite{logg2012dolfin} to compute the numerical solutions of PDEs, the Python library \textit{PyDMD} \cite{demo2018pydmd} to complete the calculations of Exact DMD,  and the Python library \textit{PyTroch} \cite{paszke2019pytorch} to train the neural networks, respectively. Besides, we employ the publicly available implementation of EDMD~\footnote{\href{https://github.com/MLDS-NUS/KoopmanDL}{https://github.com/MLDS-NUS/KoopmanDL}} , whose observable functions are  $[1,\mathbf{x},RBF(\mathbf{x})]^T$. Here, $RBF(\cdot)$ represents radial basis functions (RBF) dictionary, which consists of thin-plate RBF functions with centers placed on the training data using k-means clustering. The Xavier normal initialization scheme \cite{glorot2010understanding} is utilized to initialize the weights of all neural networks, while the biases of all nodes are set to zero. All the networks are trained by the Adam optimizer \cite{DBLP:journals/corr/KingmaB14} with an initial learning rate of $10^{-3}$. In order to find the optimal parameters of the  network, we use \textit{ReduceLROnPlateau} \cite{DBLP:journals/corr/Ruder16} to adjust the learning rate during the training process for all numerical examples.
For fairness, all the methods share the same training strategies. Denote $x$ as the ``true'' value of the states and $\hat{x}$ as its reconstruction. We use three metrics to evaluate different methods synthetically, \textit{i.e.}, the relative $L_2$ error  
\[
\text{RL2E}(t) = \frac{||\hat{x}_t-x_t||_2}{||x_t||_2},
\]
the mean squared error
\[
\text{MSE}(t) =\frac{||\hat{x}_t-x_t||^2_2}{m},
\]
and the total relative $L_2$ error
\[
\text{TRL2E} = \sqrt{\frac{\sum_{t=1}^T||\hat{x}_t-x_t||^2_2}{\sum_{i=1}^T||x_t||^2_2}}.
\]
\subsection{Fixed-point attractor}\label{sec:fixed_point}
The fixed-point attractor example \cite{lusch2018deep} is given by
\begin{equation*}
    \left\{ \begin{aligned}
        x_{t+1,1} &=\lambda x_{t,1},\\
        x_{t+1,2} &=\mu x_{t,2} + (\lambda^2-\mu)x_{t,1}^2.
    \end{aligned} 
    \right.
\end{equation*}
The initial state is chosen randomly by $x_{0,1}\sim U(0.2,4.2) $, $x_{0,2} \sim U(0.2,4.2)$ and $\lambda=0.9 , \mu=0.5$.
We divide the data set into three parts where the ratio of training, validation, and test is $60\%, 20\%$, and $20\%$, respectively. The number of neurons of each layer for the encoder network in LIR-DMD is $2,10,10,3$ and the number of neurons of  decoder network is $3,10,10,2$. This results in 345 trainable parameters for LIR-DMD. We use three ACFs for this problem. The mappings $t$ and $s$ are parameterized by FNN with three layers and the width of each layer is 1,8,2, respectively. This results in total 102 trainable parameters in FlowDMD. For  EDMD, we choose 3 RBF functions as the RBF dictionary.

\begin{figure}[H]
    \begin{center}
        \subfloat[The trajectories of $x_1$]{
            \includegraphics[width=0.4\linewidth,trim=0 0 0 0]{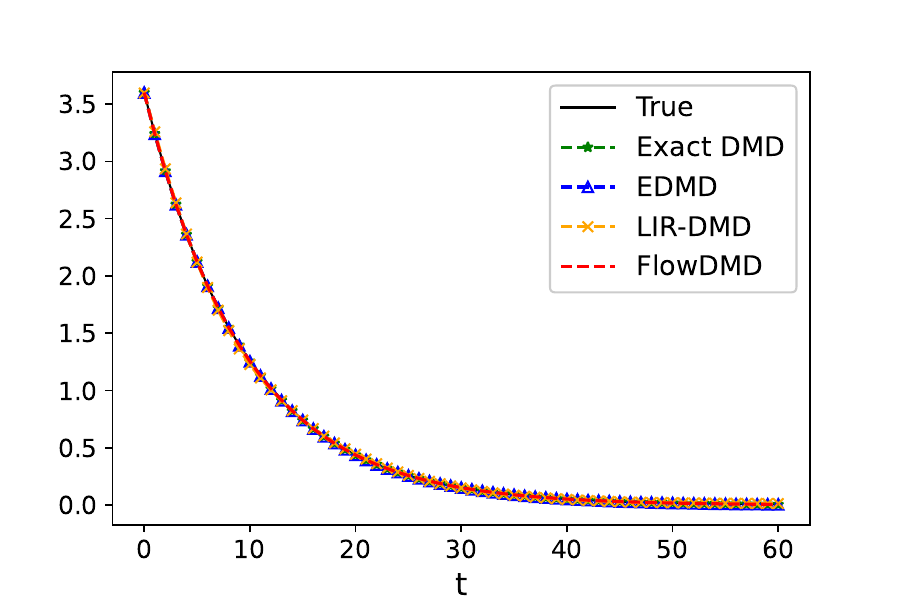}
        }\qquad
        \subfloat[The trajectories of $x_2$]{
            \includegraphics[width=0.4\linewidth,trim=0 0 0 0]{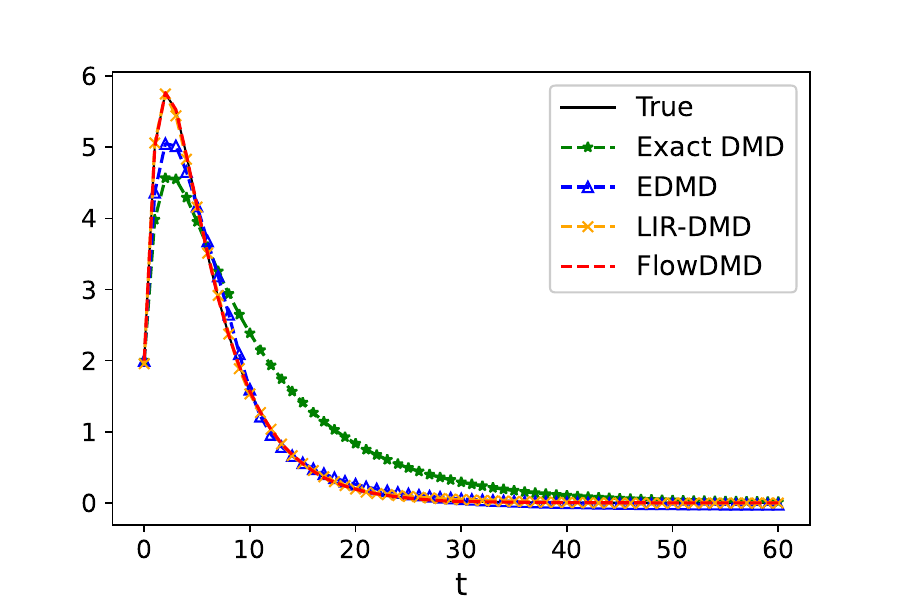}
        }\\
        \subfloat[Relative $L_2$ error]{
            \includegraphics[width=0.4\linewidth,trim=0 0 0 0]{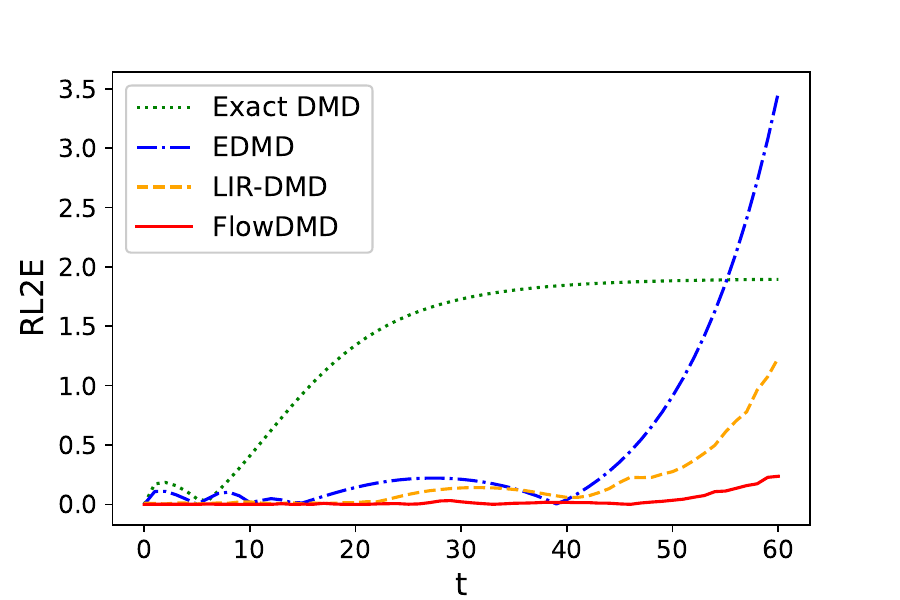}
        }\qquad
        \subfloat[Mean squared error]{
            \includegraphics[width=0.4\linewidth,trim=0 0 0 0]{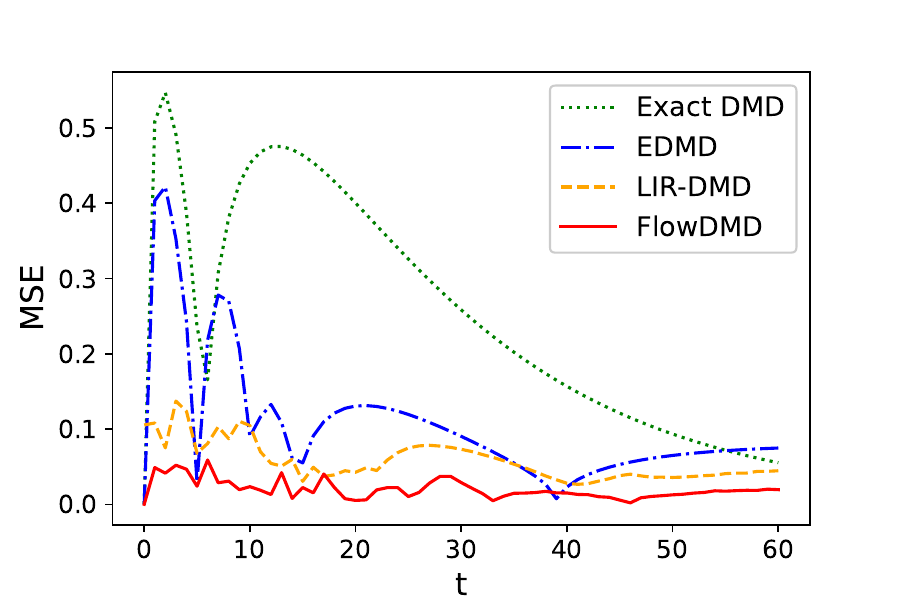}
        }
    \end{center}
    \vspace{-.4cm}
    \caption{Comparison of four methods for Example~\ref{sec:fixed_point}. The total relative $L_2$ error of the Exact DMD,  EDMD, LIR-DMD, and FlowDMD are 0.2448, 0.08,  0.0111 and 0.0018, respectively.}\label{fig:e1_1}
\end{figure}

We randomly choose one example from the test set and plot its results in Figure~\ref{fig:e1_1}. Both Figure \ref{fig:e1_1}(a) and Figure \ref{fig:e1_1}(b) show that the reconstruction calculated by LIR-DMD and FlowDMD are better than that of the Exact DMD and  EDMD. Furthermore, the difference of trajectories between LIR-DMD and FlowDMD is very small. Figure \ref{fig:e1_1}(c) and Figure \ref{fig:e1_1}(d) illustrate that the reconstruction error of FlowDMD is the smallest. In the first 30 time steps, LIR-DMD has a similar error to FlowDMD. The error of FlowDMD increases much more slowly than that of LIR-DMD for the following 30 time steps. We conclude that FlowDMD has better generalization ability than LIR-DMD.

\begin{figure}[ht]
    \begin{center}
        \subfloat{
            \includegraphics[width=0.5\linewidth,trim=0 0 0 0]{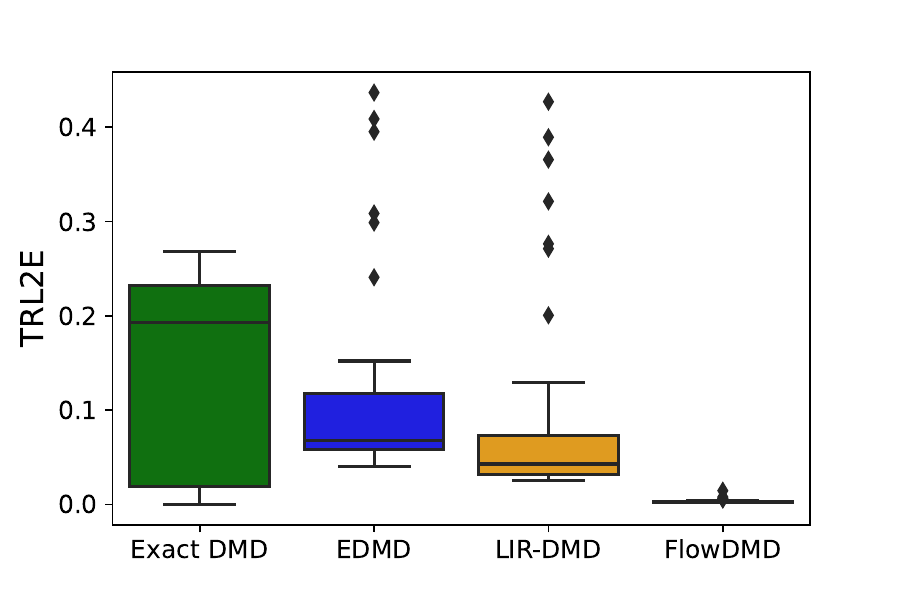}
        }
    \end{center}
    \vspace{-.6cm}
    \caption{Total relative $L_2$ error in Example~\ref{sec:fixed_point}.}\label{fig:e1_2}
\end{figure}

We use the TRL2E to evaluate the reconstruction results of trajectories for Exact DMD, EDMD, FlowDMD and LIR-DMD on 40 randomly generated examples, respectively, and the corresponding results are depicted by Figure~ \ref{fig:e1_2}. For FlowDMD, the reconstruction error is the lowest among almost all of the test examples, and the average total relative $L_2$ error is only $0.3\%$. Compared with LIR-DMD, FlowDMD has better generalization ability and learning ability of the Koopman invariant subspace.
\vspace{-.6cm}

\subsection{Burgers' equation}\label{sec:Burgers}
Consider the 1-D Burgers' equation~\cite{raissi2019physics} given by
\begin{equation*}
    \left\{ \begin{aligned}
        &\frac{\partial u}{\partial t}+u\frac{\partial u}{\partial x} = \frac{0.01}{\pi}\frac{\partial^2u}{\partial x^2}\quad x\in(-1,1), t\in(0,1],\\
        &u(1,t)=u(-1,t)=0,\\
        &u(x,0) = -\xi * sin(\pi x),
    \end{aligned} 
    \right.
\end{equation*}
where $\xi$ is a random variable that satisfies a uniform distribution $U(0.2,1.2)$. We use the finite element method with 30 equidistant grid points for the spatial discretization and the implicit Euler method with a step size of $0.01$ for temporal discretization. We generate 100 samples of $\xi$ for the initial state and compute the corresponding solutions. The examples are then divided into three parts, with proportions $60\%$ for training, $20\%$ for validation, and $20\%$ for test.
We test the performance of the Exact DMD, LIR-DMD, and FlowDMD. The rank of Exact DMD is 3 and the same rank is also used in LIR-DMD and FlowDMD to embed the Koopman linearity. The structure of the encoder network for LIR-DMD is $[30,40,50,40]$, and the decoder network is $[40,50,40,30]$ where the numbers in the brackets represent the width of each layer and we use RCFs to replace ACFs. This results in an invertible neural network of depth of 3 with one RCF block and one RCF. In each RCF, the width of each layer in FNN to parameterize the mapping $t$ is 15, 40, 15, which results in 7530 parameters in FlowDMD, whereas LIR-DMD has 10650 parameters. For  EDMD, we choose 30 RBF functions as the RBF dictionary.

\begin{figure}[H]
    \begin{center}
            \includegraphics[width=0.7\linewidth,trim=0 0 0 0]{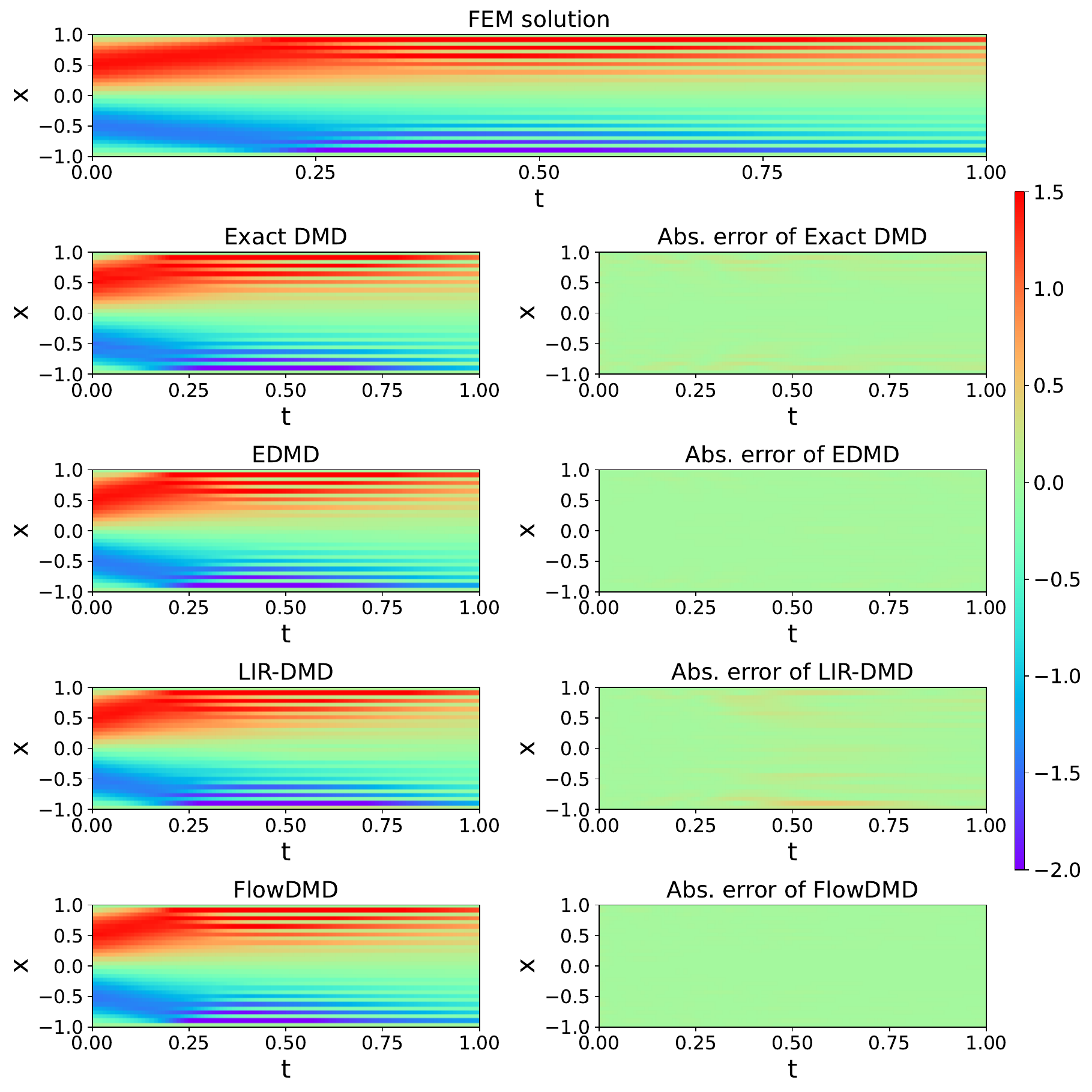}
    \end{center}
    \vspace{-.4cm}
    \caption{Comparison of four methods in Example~\ref{sec:Burgers}. The total relative $L_2$ errors for exact DMD, EDMD, LIR-DMD, and FlowDMD are 0.08,0.026, 0.119, and 0.017, respectively.}\label{fig:B_example}
\end{figure}

\begin{figure}[H]
    \begin{center}
        \subfloat[Relative $L_2$ error]{
            \includegraphics[width=0.5\linewidth]{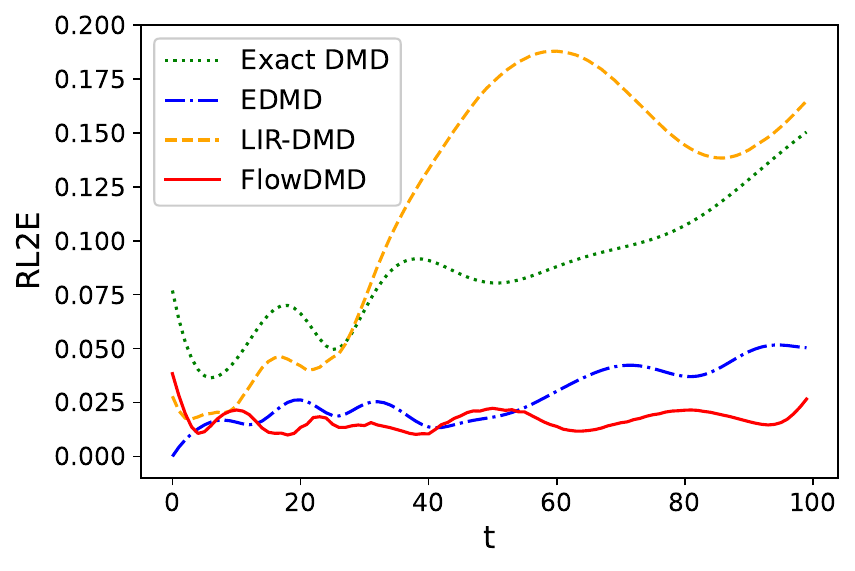}
        }
        \subfloat[Mean squared error]{
            \includegraphics[width=0.5\linewidth]{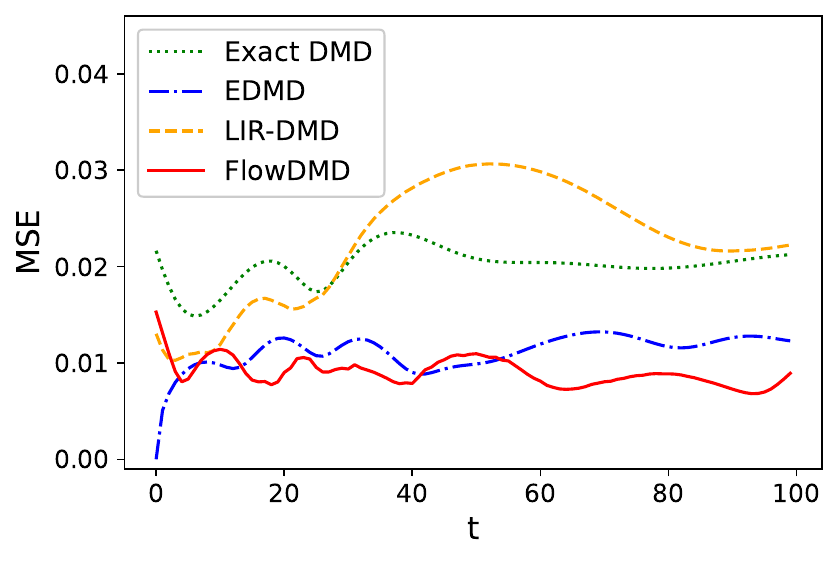}
        }
    \end{center}
    \vspace{-.6cm}
    \caption{Error of four methods for Example \ref{sec:Burgers}. }\label{fig:B_error}
\end{figure}

\begin{figure}[H]
    \begin{center}
        \subfloat{
            \includegraphics[width=0.5\linewidth]{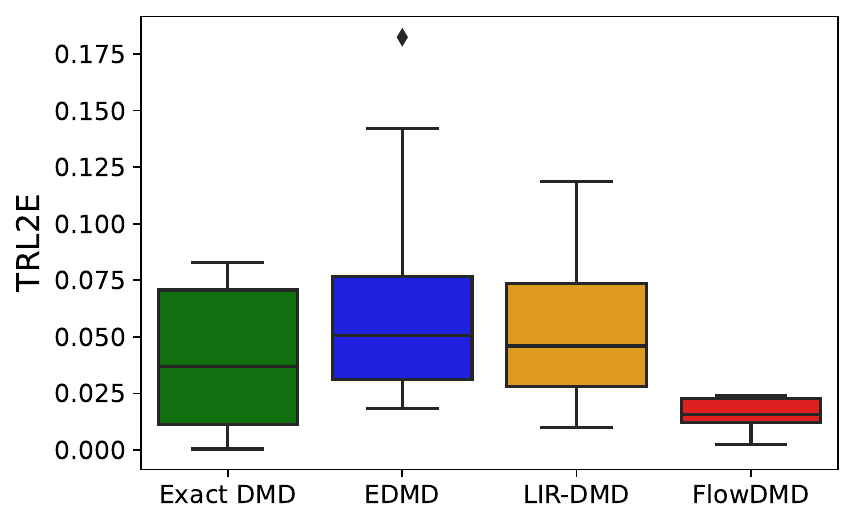}
        }
    \end{center}
    \vspace{-.6cm}
    \caption{Total relative $L_2$ error in Example~\ref{sec:Burgers}.}\label{fig:B_testdata}
\end{figure}

Figure~\ref{fig:B_example} depicts that FlowDMD has the smallest absolute reconstruction error and TRL2E. Figure~\ref{fig:B_error}(a) and Figure~\ref{fig:B_error}(b) show that the reconstruction error of Exact DMD, EDMD and LIR-DMD all increase with time, but FlowDMD maintains in a very low level.
Figure~\ref{fig:B_testdata} summarizes the TRL2E of reconstruction on all test examples and depicts that the FlowDMD has the smallest error on almost all test examples, where the average TRL2E of FlowDMD is $ 1.5\%$. For some test examples, Exact DMD has the same TRL2E with FlowDMD, but for most test examples, FlowDMD performs better than Exact DMD. The TRL2E of LIR-DMD are bigger than FlowDMD over all the test examples and are slightly better than Exact DMD for some test examples.
\vspace{-.5cm}
\vspace{-.5cm}

\subsection{Allen-Cahn equation}\label{sec:AC}
\begin{figure}[H]
        \begin{center}
                \includegraphics[width=0.75\linewidth]{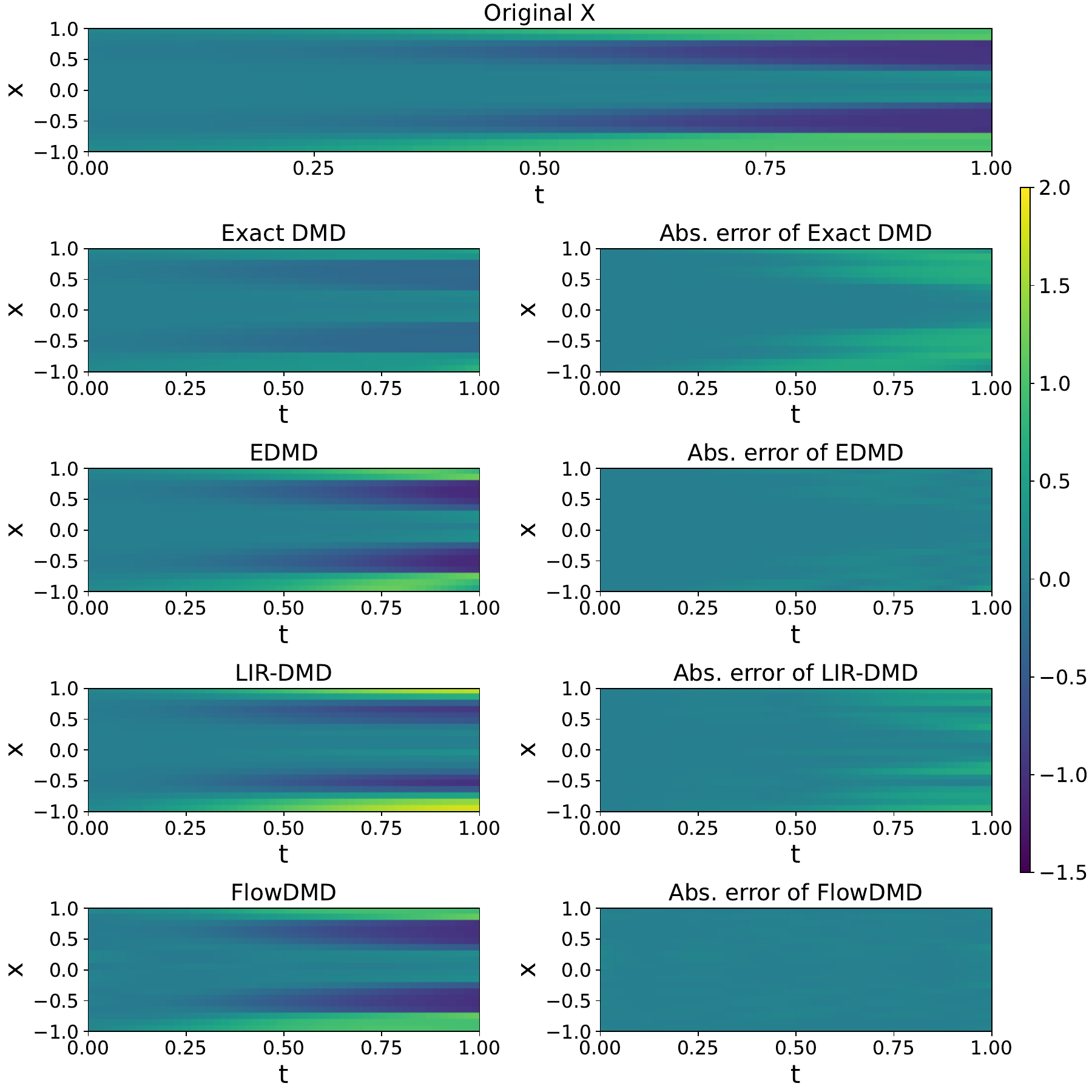}
        \end{center}
        \vspace{-.5cm}
        \caption{Comparison of four methods in Example \ref{sec:AC}. The total relative $L_2$ error for exact DMD, EDMD, LIR-DMD, and FlowDMD are 0.6129, 0.129, 0.4038, and 0.0725, respectively.}\label{fig:AC_example}
    \end{figure}

\begin{figure}[H]
        \begin{center}
            \subfloat[Relative $L_2$ error]{
                \includegraphics[width=0.45\linewidth,trim=0 0 10 0]{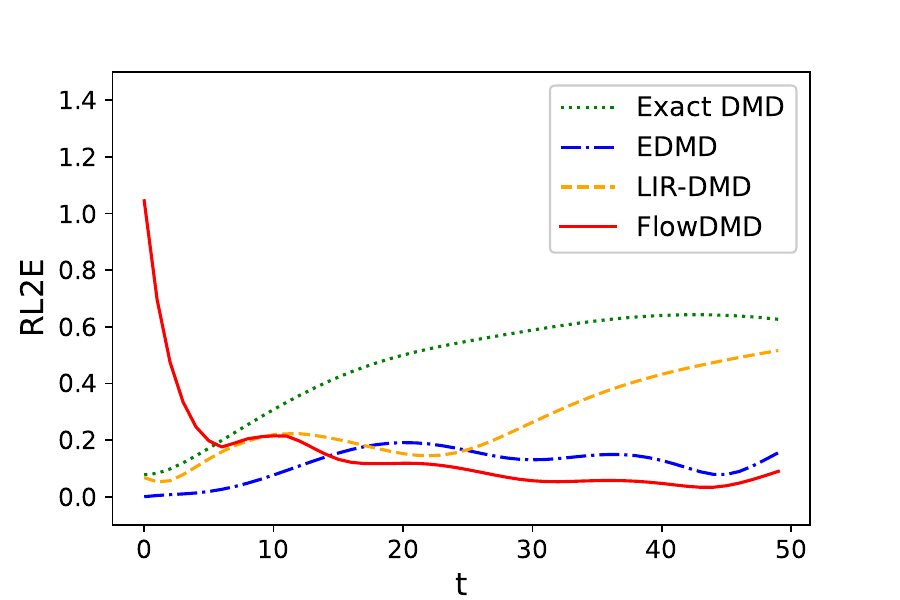}
            }\qquad
            \subfloat[Mean squared error]{
                \includegraphics[width=0.45\linewidth,trim=0 0 0 0]{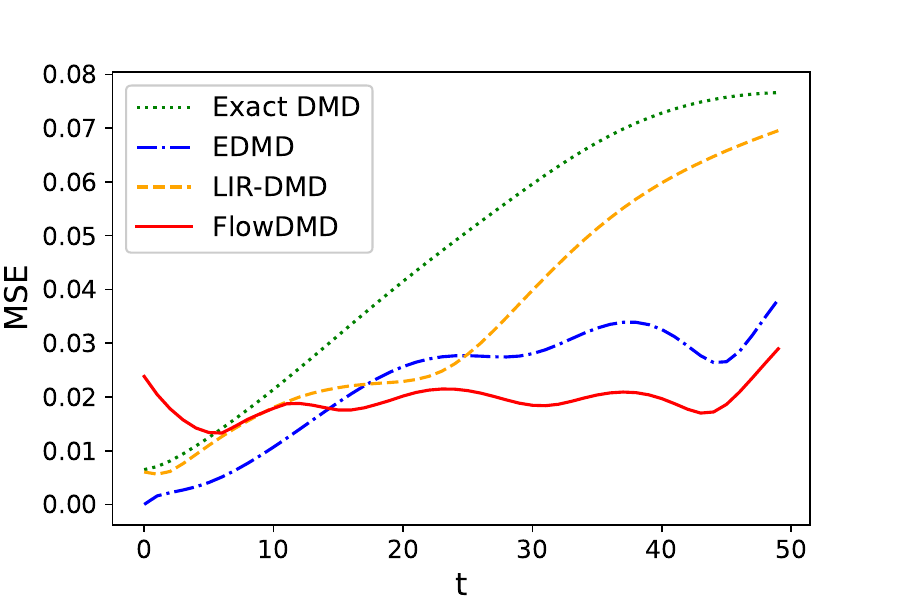}
            }
        \end{center}
        \vspace{-.6cm}
                \caption{Error of four methods for Example \ref{sec:AC}.}\label{fig:AC_error}
    \end{figure}

\begin{figure}[H]
        \begin{center}
            \subfloat{
                \includegraphics[width=0.5\linewidth]{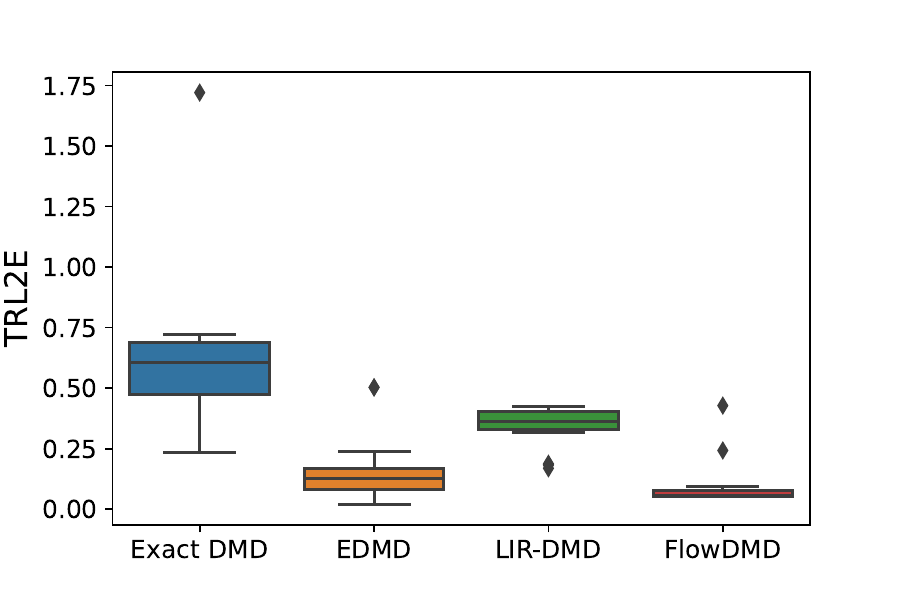}
            }
        \end{center}
        \vspace{-.6cm}
        \caption{Total relative  $L_2$ error in Example \ref{sec:AC}.}\label{fig:AC_testdata}
    \end{figure}

Consider the 1-D Allen-Cahn equation~\cite{raissi2019physics} given by
\begin{equation*}
    \left\{ \begin{aligned}
        &\frac{\partial u}{\partial t}-\gamma_1 \frac{\partial^2 u}{\partial x^2}+\gamma_2\left(u^3-u\right)=0, x \in(-1,1), t \in(0,1] ,\\
        &u(0, x)=\xi*x^2 \cos (2 \pi x) ,\\
        &u(t,-1)=u(t, 1) ,
    \end{aligned}  
    \right.
\end{equation*}
where  $\gamma_1=0.0001$, $\gamma_2=5$, and $\xi\sim \mathcal{N}(-0.1,0.04)$. We use the finite element method with 20 equidistant grid points for the spatial discretization and the implicit Euler with a step size of $0.02$ for the temporal discretization. Furthermore, we generate 100 samples of $\xi$ and use \textit{FEniCS}  to compute the numerical solutions.
The data set is segmented according to a ratio of $60\%$, $20\%$, $20\%$, respectively to be used as the training set, the validation set, and the test set. The structure of the encoder network for LIR-DMD is $[20,30,40,30 ]$ and the decoder network is $[30,40,30,20 ]$, where the numbers in the bracket indicate the width of each layer. This results in 6190 parameters for LIR-DMD. For FlowDMD, we also use RCFs to replace the ACFs.  The neural network for FlowDMD consists of one RCF block and one RCF, which results in a network with depth $L=3$. In each RCF, the width of each layer of the FNN to parameterize $t$ is 10, 20, 10. Finally, we obtain  2580 parameters for FlowDMD. The rank of Exact DMD is 3, and the same rank is also used in LIR-DMD and FlowDMD to embed the Koopman linearity. We choose 4 RBF functions as the RBF dictionary for EDMD. Results are reported in Figure~\ref{fig:AC_example}--\ref{fig:AC_testdata}.
\vspace{-.5cm}
Figure~\ref{fig:AC_example} clearly shows that FlowDMD can reconstruct the original state most accurately. It reveals that the absolute error of exact DMD, EDMD and LIR-DMD increase over time, but FlowDMD can maintain the error in a low level all the time.
In addition, numerical results show that FlowDMD is more robust and  generalizes better than Exact DMD, EDMD and LIR-DMD. Specifically, the error of the state reconstruction for four methods are given in Figure~\ref{fig:AC_error}. At the beginning time, FlowDMD has the biggest relative error  because the norm of the true state variables is too small, which leads to a large relative error. As time evolves, the error of FlowDMD reaches the lowest level among all four methods. 
\vspace{-.7cm}
    \vspace{-.2cm}
In Figure~\ref{fig:AC_testdata}, we use the test data set to evaluate the  generalization ability. The FlowDMD has almost  the smallest  TRL2E in most examples and the average of the total relative $L_2$ error is $ 9\%$. It also shows that  the fluctuation of error for FlowDMD is smaller than that of  LIR-DMD, which demonstrates that FlowDMD has a better generalization ability and is more robust than LIR-DMD.
\vspace{-.5cm}

\subsection{Sensitivity study}
Here, we study the sensitivity of FlowDMD systematically using the Allen-Cahn equation in Section~\ref{sec:AC} with respect to the following five aspects,
\begin{enumerate}
    \item The neural network initialization.
    \item The hyperparameter $\alpha$ in the loss function.
    \item The structure of neural networks.
    \item The rank $r$ used by DMD in Algorithm \ref{alg:ExactDMD}.
    \item The division index $q$ in Definition~\ref{def:CF} for CF-INN.
\end{enumerate}

\subsubsection{Sensitivity with respect to the neural network initialization}
\begin{figure}[H]
    \begin{center}
        \subfloat{
            \includegraphics[width=0.45\linewidth]{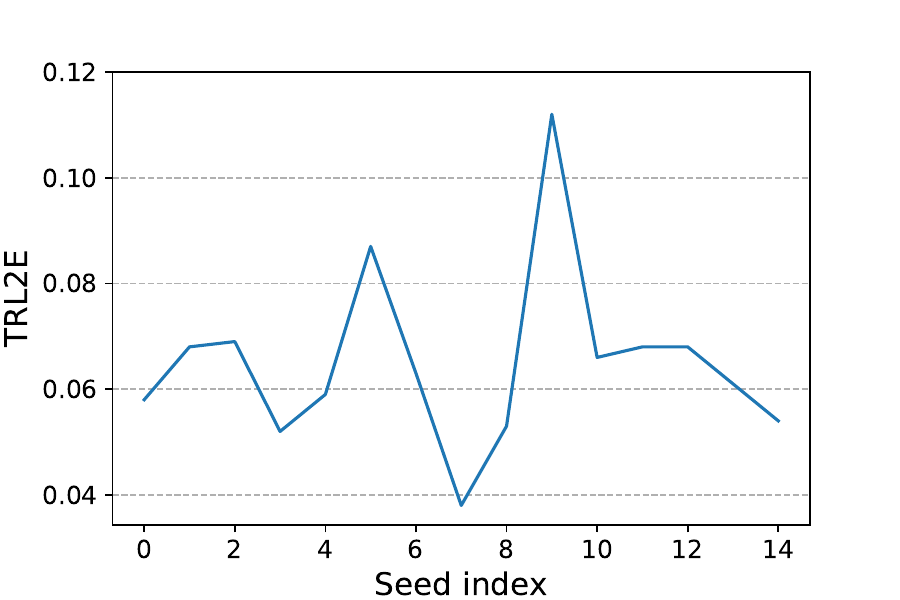}
        }
    \end{center}
    \vspace{-.5cm}
    \caption{Total relative $L_2$ error for different neural network initialization.}\label{fig:different_initialization}
\end{figure}

In order to quantify the sensitivity of FlowDMD with respect to the initialization, we consider the same data set with Section~\ref{sec:AC}. Simultaneously, we fix the structure for FlowDMD to include only one RCF block and one RCF. Each RCF has a FNN to parameterize $t$ where the width of each layer is $10,20,10$. Moreover, all FNNs use the rectified linear unit as activation functions. We use 15 random seeds to initialize models and train all the models with the same setting. In Figure~\ref{fig:different_initialization}, we report the TRL2E between the reconstructed states and the ``true" states. Evidently, the TRL2E remains stable for different initializations of neural networks, as demonstrated by the consistent results obtained within the following interval,
\[
[\mu_{TRL2E}-\sigma_{TRL2E},\mu_{TRL2E}+\sigma_{TRL2E}]=[6.5\times10^{-2}-1.6\times10^{-2},6.5\times10^{-2}+1.6\times10^{-2}]
\]
\vspace{-.5cm}

\subsubsection{Sensitivity with respect to $\alpha$}
We utilize the same training set with Section~\ref{sec:AC} and select $\alpha$ from the list $[0.01,0.1,1,10,100]$. As shown in Table \ref{tab:different_alpha}, the different weights $\alpha$ in the loss function have little influence on the final results. We observe that the error is minimized when $\alpha=10$, which suggests the use of an adaptive weight selection algorithm. The gradient flow provided by the neural tangent kernel \cite{WANG2022110768} can be employed to adjust the weight $\alpha$ and accelerate the training process, and we leave this for our future work.
\begin{table}[ht]
    \caption{Total relative $L_2$ error for different $\alpha$.}\label{tab:different_alpha}
    \vspace{-.5cm}
    \begin{center}
        \begin{tabular}{cccccc}
            \hline  
            $\alpha$ & 0.01 &0.1& 1& 10&100\\
            \hline
            TRL2E &6.2e-02& 6.8e-02& 8.2e-02& 3.2e-02& 6.9e-02\\            
            \hline  
            \end{tabular}
    \end{center}
\end{table}
\vspace{-.6cm}
\subsubsection{Sensitivity with respect to the structure of neural networks}
We study the impact of the number of RCFs and the number of neurons in the FNN to parameterize the mapping $t$ on the performance of the FlowDMD. Specifically, the sensitivity of FlowDMD is being quantified with respect to two parameters: the number of RCFs ($N_f$) and the number of neurons ($N_n$) in the middle layer of the FNN. Here, the FNN used to parameterize $t$ is restricted to a three layer structure of $[10, N_n, 10]$. The results are summarized in Table \ref{tab:different_structure}, which indicate that the reconstruction of FlowDMD has little to do with its structure while adding more neurons or more RCFs will not improve the final results to a big extent.
\begin{table}[ht]
    \caption{Total relative $L_2$ error for different structures of networks in FlowDMD.}\label{tab:different_structure}
    \vspace{-.5cm}
    \begin{center}
        \begin{tabular}{cccccc}
            \hline 
            \diagbox{$N_n$}{$N_f$} & 2 &3& 4& 8&12\\
            \hline
            10 &5.6e-02& 7.7e-02& 7.1e-02& 8.4e-02& 5.4e-02\\            
            20 &6.4e-02& 7.6e-02& 6.4e-02& 7.9e-02& 8.6e-02\\            
            40 &7.1e-02& 8.4e-02& 4.1e-02& 4.8e-02& 10.3e-02\\            
            80 &4.0e-02& 7.6e-02& 7.7e-02& 8.0e-02& 6.6e-02\\            
            160 &8.9e-02& 4.2e-02& 8.0e-02& 8.3e-02& 6.9e-02 \\
            \hline  
            \end{tabular}
    \end{center}
\end{table}
\vspace{-.6cm}
\subsubsection{Sensitivity with respect to the rank of DMD}
As we increase the rank $r$ used for the DMD computations in Algorithm~\ref{alg:ExactDMD}, we include more information, but the computation time also increases.  In this study, we investigate how the DMD rank affects the model and its reconstruction. The results in Table \ref{tab:different_r} show that as we increase the rank $r$, the corresponding error decreases rapidly.
\begin{table}[htbp]
    \caption{Total relative $L_2$ error for different low rank dimension in Algorithm \ref{alg:ExactDMD}.}\label{tab:different_r}
    \vspace{-.5cm}
    \begin{center}
        \begin{tabular}{cccccc}
            \hline
            $r$ & 1 &3& 5& 7&9\\
            \hline
            TRL2E &17.4e-02& 6.8e-02& 6.7e-02& 9e-03& 3e-03\\            
            \hline
            \end{tabular}
    \end{center}
\end{table}
\vspace{-.6cm}

\subsubsection{Sensitivity with respect to $q$}
\begin{figure}[H]
    \begin{center}
        \subfloat{
            \includegraphics[width=0.4\linewidth]{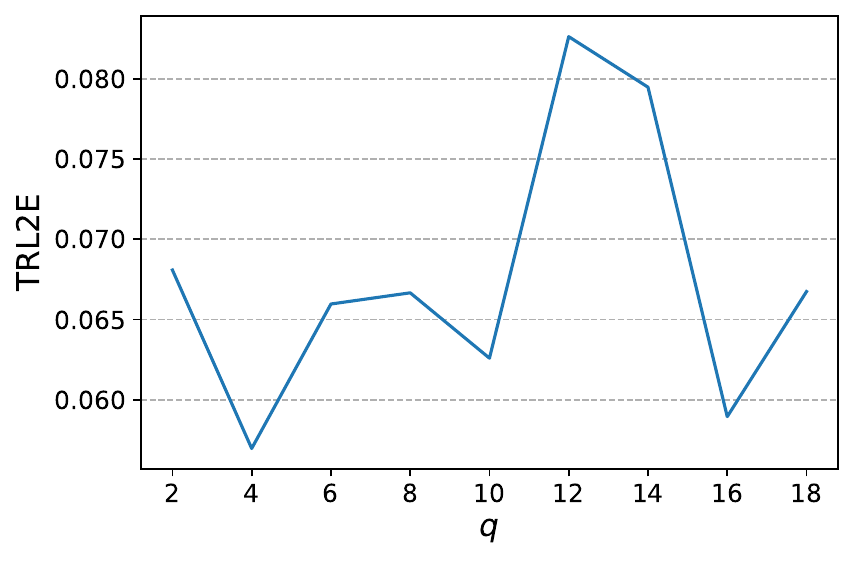}
        }
    \end{center}
    \vspace{-.5cm}
    \caption{Total relative $L_2$ error for different division indices.}\label{fig:different_q}
\end{figure}

To investigate how the division index $q$ affects the performance of FlowDMD, we report the mean value of TRL2E for different values of $q$ in Figure~\ref{fig:different_q}. As Figure~\ref{fig:different_q} shows, the setting of $q=4$ minimizes the TRL2E and $q=12$ results in the largest TRL2E. However, the total relative $L_2$ error always remains in a low and stable level by FlowDMD.

\section{Conclusion}
\label{sec-Conclusion}
In this paper, we introduced the FlowDMD framework to learn both the observable functions and reconstruction functions for the Koopman embedding, which is implemented through coupling flow invertible neural network. Our method gives more accurate approximations of the Koopman operator than state-of-the-art methods. Our FlowDMD is structurally invertible, which simplifies the loss function and improves the accuracy of the state reconstruction. Numerical experiments show that our approach is more accurate, efficient, and interpretable than the state-of-the-art methods.



\section*{Acknowledgments}
We would like to thank Mengnan Li and Lijian Jiang for sharing their code. We also would like to thank the anonymous referees that help to improve the quality of this paper.




\bibliographystyle{elsarticle-num-names} 
\bibliography{ref, ref_myh}
\end{document}